\documentclass{article}

\newcommand{\rme}{\mathrm{e}}

\newcommand{\dx}{\Delta x}

\newcommand{\wh}{\widehat}

\newcommand{\Da}{\mathrm{Da}}

\newcommand{\ra}{\mathrm{a}}
\newcommand{\rd}{\mathrm{d}}

\newcommand{\ru}{\mathrm{u}}

\newcommand{\p}{\partial}
\newcommand{\wt}{\widetilde}

\usepackage{amsmath,amsfonts,amssymb,enumerate, tikz, graphicx,float,makecell,subfig, soul, url, hyperref,adjustbox}
\usetikzlibrary{arrows,decorations.markings}

\numberwithin{equation}{section}
\numberwithin{figure}{section}
\numberwithin{table}{section}

\begin{document}

\title{Diffusion-Limited Reactions in Nanoscale Electronics \thanks{Submitted to the editors October 18, 2017.} }

\author{Ryan M. Evans 
\thanks{Applied and Computational Mathematics Division, National Institute of Standards and Technology, Gaithersburg, MD 20899, USA. (\href{ryan.evans@nist.gov}{ryan.evans@nist.gov}, \href{anthony.kearsley@nist.gov}{anthony.kearsley@nist.gov}) \textbf{Funding:} The first author was supported by the National Research Council through a postdoctoral fellowship.
}
\and Arvind Balijepalli 
\thanks{Engineering Physics Division, National Institute of Standards and Technology, Gaithersburg, MD 20899, USA. (\href{arvind.balijepalli@nist.gov}{arvind.balijepalli@nist.gov})}%
\and Anthony J. Kearsley\footnotemark[2]
}

\maketitle

\begin{abstract}
A partial differential equation (PDE) was developed to describe time-dependent ligand-receptor interactions for applications in biosensing using field effect transistors (FET). The model describes biochemical interactions at the sensor surface (or biochemical gate) located at the bottom of a solution-well, which result in a time-dependent change in the FET conductance. It was shown that one can exploit the disparate length scales of the solution-well and biochemical gate to reduce the coupled PDE model to a single nonlinear integrodifferential equation (IDE) that describes the concentration of reacting species.  Although this equation has a convolution integral with a singular kernel, a numerical approximation was constructed by applying the method of lines.  The need for specialized quadrature techniques was obviated and numerical evidence strongly suggests that this method achieves first-order accuracy.  Results reveal a depletion region on the biochemical gate, which non-uniformly alters the surface potential of the semiconductor.


\end{abstract}


\section{Introduction} \label{Section: Introduction}

The ability to tailor therapies to individuals or specific subsets of a population to deliver personalized care has the potential to fundamentally remake healthcare delivery. The most promising therapeutic candidates for such targeted care are new classes of biologic drugs based on naturally occurring molecules, made possible due to rapid advances in genomics and proteomics \cite{Fosgerau:2015kt, Walsh:2014kw}. Importantly, such therapies can be safer and yield better outcomes at lower doses when treating debilitating conditions such as diabetes, Alzheimer\textcolor[rgb]{0,0,1}{'}s disease, or certain cancers \cite{ Binukumar:2014hl,Dhavan:2001gd}. The widespread use of personalized care is currently limited by our ability to routinely measure pathology in individuals including biomarkers, metabolites, tissue histology, and gene expression. Moreover, existing clinical diagnostics are cumbersome, require specialized facilities, can take days to weeks to perform, and are in many cases prohibitively expensive. This has led to the development of new portable detection tools including antibody-based lateral flow assays \cite{Drain:2016fs, Lawn:2016ed}, microelectromechanical sensor (MEMS) based resonators that can detect binding of biomarkers to the sensor surface \cite{ Ilic:2004go, Johannsmann:2015kq,  Naik:2009gu, Rodahl:1995ck}, surface plasmon resonance \cite{ Knoll:1998df,Wang:2010bo}, ring cavity resonators \cite{Armani:2003gx, Boyd:2001tg, Su:2015fo}, and electronic measurements with field effect transistors (FET) \cite{ Cui:2001iy, Mohanty:2014wy, Pouthas:2004cj, Wang:2005ij}. The latter are particularly well-suited for biomarker measurements due their high charge sensitivity and direct signal transduction, allowing label-free measurements at physiological concentrations. Furthermore, by leveraging semi-conductor processing techniques, measurements with FETs can be made massively parallel, cost-effective, and portable.

\begin{figure}[tbhp!]
\centering
\begin{tikzpicture}
 
\draw [thick] (-3,0)--(-3,3);
\draw [thick] (-3,0)--(-.5,0);
\draw [ultra thick,gray] (-.5,0)--(.5,0);
\draw [thick] (.5,0)--(3,0);
\draw [thick] (3,0)--(3,3);

\draw [ultra thick] (-.5,-.05)--(.5,-.05);


\draw [thick] (-.75,0)--(-.75,.-.25);
\draw [thick] (-.75,-.25)--(-.25,-.25);
\draw [thick] (-.25, -.25)--(-.25,-.07);

\node [left] (1) at (-1.5,-.125) {\tiny Source};
\node (2) at (-.75,-.125) {};
\path[->,>=stealth',shorten >=1pt] 
(1) edge  (2);


\draw [thick] (.25,-.07)--(.25,-.25);
\draw [thick] (.25,-.25)--(.75,-.25);
\draw [thick] (.75,-.25)--(.75,0);


\node [right] (3) at (1.5,-.125) {\tiny Drain};
\node (4) at (.75,-.125) {};
\path[->,>=stealth',shorten >=1pt] 
(3) edge  (4);


\node (7) at (.25,-.05) {};
\node [right] (8) at (1.25,-.75) {\tiny Semiconductor channel};
\path[->,>=stealth',shorten >=1pt]
(8) edge  (7);

\node (9) at (-.3,.02) {};
\node [  left] (10) at (-1,.6) {\tiny Biochemical gate};
\path[->,>=stealth',shorten >=1pt]
(10) edge  (9);


\draw (-.25,0)--(-.25,.2);
\draw (-.25,.2)--(-.4,.3);
\draw (-.25,.2)--(-.1,.3);

\draw (.25,0)--(.25,.2);
\draw (.25,.2)--(.4,.3);
\draw (.25,.2)--(.1,.3);


\draw (-.3,2.65) circle [radius=.1];
\draw (.1,2.35) circle [radius=.1];
\draw (.27,2.75) circle [radius=.1];
\draw (.5,2.3) circle [radius=.1];
\draw (-.4,2.21) circle [radius=.1];
\draw (-.1,1.95) circle [radius=.1];
\draw (.1,1.55) circle [radius=.1];
\draw (.35,1.95) circle [radius=.1];
\draw (.375,1.35) circle [radius=.1];
\draw (-.35,1.65) circle [radius=.1];
\draw (0,1.15) circle [radius=.1];
\draw (-.4,1.25) circle [radius=.1];
\draw (.35,.9) circle [radius=.1];
\draw (-.275,.825) circle [radius=.1];
\draw (.125,.6) circle [radius=.1];
\draw (-.35,.4725) circle [radius=.1];

\draw (-.75,2.7) circle [radius=.1];
\draw (-.85,2.3) circle [radius=.1];
\draw (-1.2,2.5) circle [radius=.1];
\draw (-1.6,2.65) circle [radius=.1];
\draw (-1.6,2.25) circle [radius=.1];
\draw (-2,2.5) circle [radius=.1];
\draw (-2.2,2.2) circle [radius=.1];
\draw (-2.5,2.6) circle [radius=.1];
\draw (-2.6,2.2) circle [radius=.1];

\draw (.75,2.7) circle [radius=.1];
\draw (.85,2.3) circle [radius=.1];
\draw (1.2,2.5) circle [radius=.1];
\draw (1.6,2.65) circle [radius=.1];
\draw (1.6,2.25) circle [radius=.1];
\draw (2,2.5) circle [radius=.1];
\draw (2.2,2.2) circle [radius=.1];
\draw (2.5,2.6) circle [radius=.1];
\draw (2.6,2.2) circle [radius=.1];

\draw (-.75, 1.9) circle [radius=.1];
\draw (-.8, 1.55) circle [radius=.1];
\draw (-1.2,1.6) circle [radius=.1];
\draw (-1.15, 2.05) circle [radius=.1];
\draw (-1.55,1.9) circle [radius=.1];
\draw (-1.6,1.5) circle [radius=.1];
\draw (-1.9, 2) circle [radius=.1];
\draw (-2, 1.65) circle [radius=.1];
\draw (-2.4, 1.8) circle [radius=.1];
\draw (-2.4, 1.4) circle [radius=.1];
\draw (-2.7,1.6) circle [radius=.1];

\draw (.75, 1.9) circle [radius=.1];
\draw (.8, 1.55) circle [radius=.1];
\draw (1.2,1.6) circle [radius=.1];
\draw (1.15, 2.05) circle [radius=.1];
\draw (1.55,1.9) circle [radius=.1];
\draw (1.6,1.5) circle [radius=.1];
\draw (1.9, 2) circle [radius=.1];
\draw (2, 1.65) circle [radius=.1];
\draw (2.4, 1.8) circle [radius=.1];
\draw (2.4, 1.4) circle [radius=.1];
\draw (2.7,1.6) circle [radius=.1];

\draw (-.75,1.2) circle [radius=.1];
\draw (-.65,.8) circle [radius=.1];
\draw (-1.05,1) circle [radius=.1];
\draw (-1.15,1.3) circle [radius=.1];
\draw (-1.5,1.15) circle [radius=.1];
\draw (-2,1.35) circle [radius=.1];
\draw (-1.85,1.1) circle [radius=.1];
\draw (-1.65,.85) circle [radius=.1];
\draw (-2.1,.85) circle [radius=.1];
\draw (-2.3,1.1) circle [radius=.1];
\draw (-2.6,1) circle [radius=.1];

\draw (.75,1.2) circle [radius=.1];
\draw (.65,.8) circle [radius=.1];
\draw (1.05,1) circle [radius=.1];
\draw (1.15,1.3) circle [radius=.1];
\draw (1.5,1.15) circle [radius=.1];
\draw (2,1.35) circle [radius=.1];
\draw (1.85,1.1) circle [radius=.1];
\draw (1.65,.85) circle [radius=.1];
\draw (2.1,.85) circle [radius=.1];
\draw (2.3,1.1) circle [radius=.1];
\draw (2.6,1) circle [radius=.1];

\draw (1.35, .8) circle [radius=.1];
\draw (1.05, .6) circle [radius=.1];
\draw (.6, .5) circle [radius=.1];
\draw (1.5, .5) circle [radius=.1];
\draw (1.9,.6) circle [radius=.1];
\draw (2.25,.6) circle [radius=.1];
\draw (2.5,.75) circle [radius=.1];
\draw (2.75,.5) circle [radius=.1];
\draw (2.5,.25) circle [radius=.1];
\draw (2,.3) circle [radius=.1];
\draw (1.5,.15) circle [radius=.1];
\draw (1.15,.25) circle [radius=.1];
\draw (.8,.2) circle [radius=.1];

\end{tikzpicture}
\caption{Schematic of biomarker measurements with a field effect transistor (FET).  Ligand molecules injected at the top of the solution-well diffuse  and bind with receptors immobilized on the FET biochemical gate.  This schematic is not drawn to scale. In particular, the  width of the solution-well is of  on the order of millimeters, and substantially larger than the size of the biochemical gate, which spans micrometers. See Table \ref{Table: parameter values} for exact parameter values.}
\label{Figure: instrument geometery}
\end{figure}
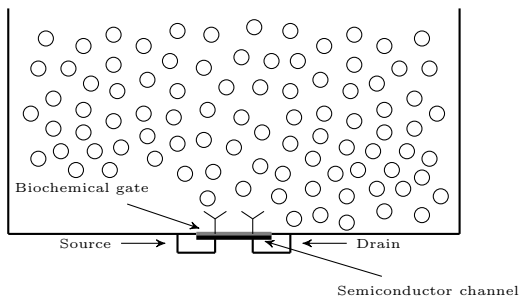

 A FET is a three-terminal device represented in Figure \ref{Figure: instrument geometery}. A semiconductor channel between the source and drain terminals conducts a current that is strongly modulated by an electrostatic potential applied to the gate. Biomarkers in aqueous solution exhibit a well-defined electrostatic surface potential \cite{Cardone:2013hh, Henrich:2010ds} arising from charged hydrophilic residues that interact with water. When these molecules adsorb to the FET biochemical gate, they strongly modulate the channel current proportional to the magnitude of their surface potential. This allows FETs to be used to detect and quantify adsorbed biomarkers in solution. Furthermore, functionalizing the FET, by attaching molecules to the gate surface that have a high inherent affinity for biomarkers of interest (see Figure \ref{Figure: instrument geometery}), allows measurements with high specificity that are tailored to one or more biomarkers of interest. 

An accurate and dynamical model of receptor ligand interactions at the biochemical gate is a critical component in maximizing the sensitivity of FET-based measurements. Specifically, quantitative descriptions of the distribution of adsorbed ligands and their surface potentials can be combined with a model of the semi-conductor physics to allow predictions of the measured signal. This in turn can be used to optimize sensor design, particularly the geometry of the biochemical gate. Of particular interest is a quantitative description of the coupling between bound ligand evolution and diffusion.  To the authors' knowledge this is a previously unexplored area of mathematical inquiry, though  Poisson-Boltzman approaches to model sensor physics have been explored.  For example in \cite{heitzinger2010multiscale} Heitzinger \textit{et al}. use the Poisson-Boltzman equation to  develop a multiple-scale model for the electric potential distribution within semiconductors of planar and nanowire field\textcolor[rgb]{0,0,1}{-}effect biosensors.  Therein, the authors model these devices using three layers: a semiconductor layer, a dielectric layer, and a discrete layer of biomolecules immobilized on the dieletric layer.  Homogenization techniques are employed to reconcile the biomolecule length scale with the semiconductor length scale, and interface conditions for the biomolecule-dielectric interface are derived. It must be noted that there are several important differences between \cite{heitzinger2010multiscale} and the present manuscript.  Perhaps the most important is that while \cite{heitzinger2010multiscale} focuses on the electric potential distribution within the semiconductor channel, the present manuscript models the coupling between reaction and diffusion.  Furthermore, while the authors of \cite{heitzinger2010multiscale} model the biomolecule layer with a discrete number of biomolecules and use homogenization techniques, in the present manuscript a continuum perspective is presented.  Finally,  \cite{heitzinger2010multiscale} assumes a \textit{steady} distribution of biomolecules on the dielectric layer, while the present manuscript concerns the time-evolution of $B(x,t)$, which is experimentally measurable.

The authors of  and  \cite{landheer2005model} employ the one-dimensional Poisson-Boltzman equation to model the electrostatic potential from a layer of biological macromolecules on the biochemical gate of a metal-oxide-semiconductor transistor. In contrast, \cite{baumgartner2011analysis} uses a three-dimensional model of the electric potential in semiconductor channel, and couples the aqueous and semiconductor regions through interface conditions obtained from Monte-Carlo simulations, which provide an estimate of the charge distribution from adsorbed biomolecules on the biochemical gate.

In \cite{tulzer2013kinetic} Heitzinger, Mauser, and Ringhofer calculate numerical values for the kinetic parameters governing adsorption and desorption processes of $\mathrm{CO}$ at a $\mathrm{SnO}{}_2$ single-nanowire gas sensor.  The authors adopt a continuum perspective by modeling surface reactions on a single-nanowire gas sensor through a set of differential equations.  However, in \cite{tulzer2013kinetic} the authors simply apply the well-stirred kinetics approximation in which gaseous carbon monoxide transport is completely divorced from adsorption and desorption processes at the surface.  This reduces their model to a set of nonlinear ordinary differential equations (ODE), which can be used to estimate kinetic rate constants involved in the reaction of interest.

In the present manuscript a quantitative description of the coupling between reaction and diffusion in FETs is developed.  In particular, we consider the experimentally relevant limit of very low ligand concentrations---\textit{i.e.}, on the order of pico- to femtomolar concentrations---and very fast assocation rates.   This problem is particularly challenging due to the disparate time and length scales involved. For example, the length scales span three orders of magnitude, ranging from order of millimeters for the solution-well to micrometers for the biochemical gate.   Combining this fact with the diffusion-limited nature of the kinetics under consideration leads to the conclusion that the time-evolution of the reacting species concentration depends heavily on a diffusive boundary layer near the surface.

In Section \ref{Section: Governing Equations} a mathematical modeled is developed that describe diffusion of ligand molecules through the solution-well onto the biochemical gate.  In Subsection \ref{Subsection: Mathematical Model} the governing equations are presented, and it is shown that there are multiple time and length scales associated with the experiment.  In Subsection \ref{Subesction: Integrodifferential Equation Reduction}, complex analysis techniques are employed to reduce the coupled PDE system to a single nonlinear integrodifferential equation (IDE)  for the reacting species concentration.  A quadrature-free numerical solution based on the method of lines is developed in Section \ref{Section: Numerical Method}, where it is shown that this method achieves first-order accuracy despite the presence of a convolution integral with a singular kernel.  Results and their physical interpretations are discussed in Section \ref{Section: Results and Discussion}, and concluding remarks are given in Section \ref{Section: Conclusions}. 






\section{Governing Equations} \label{Section: Governing Equations}

\subsection{Mathematical Model}
\label{Subsection: Mathematical Model}

Consider the geometry in Figure \ref{Figure: instrument geometery}, and a rectangular domain, $(\wt{x},\wt{y})\in\ [0,\wt{L}]\times[0,\wt{H}]$, with the origin $(0,0)$ located at the lower-left corner of the well.  The parameters $\wt{L}$ and $\wt{H}$ are the height and length of the well respectively; for parameter values see Table \ref{Table: parameter values}.  Throughout the manuscript tildes are used to denote dimensional quantities.   Receptors are confined to the biochemical gate, which occupies the very narrow region $(\wt{x},\wt{y})\in [-\wt{l}_{\mathrm{s}}/2+\wt{L}/2,\ \wt{L}/2+\wt{l}_{\mathrm{s}}]\times 0=[\wt{x}_{\mathrm{min}},\wt{x}_{\mathrm{max}}]\times 0$, where $\wt{l}_{\mathrm{s}}$ denotes length of the biochemical gate and $[\wt{x}_{\mathrm{min}},\wt{x}_{\mathrm{max}}]:=[-\wt{l}_{\mathrm{s}}/2+\wt{L}/2,\ \wt{L}/2+\wt{l}_{\mathrm{s}}]$.  It is important to note that while the length scale of the well is on the order of millimeters, the length scale of the biochemical gate is on the order of micrometers.

\begin{table}[thbp!]
 \centering
 \caption{Bounds for dimensional and dimensionless parameters are given below. }
\begin{adjustbox}{angle=270}

 \begin{tabular}{llll}
 \hline 
 \multicolumn{2}{c}{Dimensional Parameters} & \multicolumn{2}{c}{Dimensionless Parameters}\\\hline
Parameter & Range &  Parameter & Range\\ \hline
$\wt{D}\ (\mathrm{cm}^2/\mathrm{s})$ & $10^{-6}$  &$D_{\mathrm{w}}$  &$2.5\times 10^{-2}$ to $2.5\times 10^2$ \\
$\wt{k}_{\ra}\ (\mathrm{cm}^3\cdot(\mathrm{mol}\cdot \mathrm{s})^{-1})$ & $ 10^{11}$ to $10^{12}$   &$D$  & $4\times 10^{3}$ to $4\times 10^7$ \\
$\wt{k}_{\rd}\ (\mathrm{s}^{-1})$ & $10^{-5}$ to $1$    &$\Da_{\mathrm{w}}$ &$1.33\times 10^3$ to $2.66\times 10^3$ \\
$\wt{C}_{\ru}\ (\mathrm{mol}\cdot\mathrm{cm}^{-3})$ &$ 10^{-18}$ to $10^{-15}$   & $\Da$ &$3.32$ to $66.42$  \\
$\wt{R}_{\mathrm{t}}\ (\mathrm{mol}\cdot\mathrm{cm}^{-2})$ &  $6.6422\times10^{-14}$ to $1.3284\times 10^{-13}$&  $K$  &$10^{-2}$ to $10^6$ \\
$\wt{H}\ (\mathrm{cm})$ &  $0.2$   &$\epsilon$   & $0.4$ \\
$\wt{L}\ (\mathrm{cm})$ &  $0.5$  & $l_{\mathrm{s}}$ & $10^{-3}$  \\ 
$\wt{l}_s\ (\mathrm{cm})$ &  $5\times 10^{-4}$ &   &   \\ 
 \end{tabular}
 \label{Table: parameter values}
 
 \end{adjustbox}
 \end{table}

Assuming that ligand molecules are continuously and uniformly injected at the top of the well, ligand transport is governed by the diffusion equation expressed in dimensionless form as:\begin{subequations}
\begin{align}
&\frac{\p C}{\p t}=D_{\mathrm{w}}\left(\epsilon^2\frac{\p^2 C}{\p \overline{x}^2}+\frac{\p^2 C}{\p \overline{y}^2}\right),\label{C diff w}\\
&C(\overline{x},\overline{y},0)=0,\label{C ic w}\\
&C(\overline{x},1,t)=1\label{inj}\\
&\frac{\p C}{\p \overline{x}}(0,\overline{y},t)=\frac{\p C}{\p \overline{x}}(1,\overline{y},t)=0.\label{C no flux w}
\end{align} \label{C w inside}\end{subequations}
Equation (\ref{C diff w}) is the diffusion equation, (\ref{C ic w}) is the initial condition, (\ref{inj}) is the uniform injection condition, and (\ref{C no flux w}) are no-flux conditions which hold on the sides of the well.  In writing (\ref{C diff w})--(\ref{C no flux w}), we have nondimensionalized the spatial variables $\wt{x}$ and $\wt{y}$ using the well dimensions by setting $\overline{x}=\wt{x}/\wt{L}$ and $\overline{y}=\wt{y}/\wt{H}$.  Additionally, since we are interested in reaction dynamics on the sensor surface, the time variable has been scaled by the forward reaction rate $t=\wt{k}_{\mathrm{a}}\wt{C}_{\mathrm{u}}\wt{t}$. Here $\wt{C}_{\mathrm{u}}$ is the uniform injection concentration at the top of the well.  In (\ref{C diff w})--(\ref{C no flux w}) $\epsilon=O(1)$ is the aspect ratio, and
\begin{equation}
D_{\mathrm{w}}=\frac{\wt{D}}{\wt{H}^2\wt{k}_{\ra}\wt{C}_{\mathrm{u}}}=\frac{\wt{D}/{\wt{H}^2}}{\wt{k}_{\ra}\wt{C}_{\mathrm{u}}}\label{Dw}
\end{equation}
is a dimensionless constant that scales the diffusive time, $\wt{D}/\wt{H}^2$, to the forward reaction time, $\wt{k}_{\mathrm{a}}\wt{C}_{\mathrm{u}}$.   The subscript $\mathrm{w}$ indicates that the independent variables are scaled with the well dimensions.   It is seen in Table \ref{Table: parameter values} that $D_{\mathrm{w}}=O(10^{-3})$ to  $O(10)$ which implies that the the reaction at the biochemical gate is diffusion-limited, as expected for femtomolar ligand concentrations $\wt{C}_{\mathrm{u}}$. 


To state the bottom boundary condition associated with (\ref{C diff w})--(\ref{C no flux w}) we observe that when $(\overline{x},\overline{y})\not\in [\overline{x}_{\mathrm{min}},\overline{x}_{\mathrm{max}}]\times 0$ there is no flux through the surface of the well, while when $(\overline{x},\overline{y})\in [\overline{x}_{\mathrm{min}},\overline{x}_{\mathrm{max}}]\times 0$ the diffusive flux normal to the binding surface is used in forming the bound ligand .  These two conditions are expressed compactly as:
\begin{equation}
(\mathbf{n}\cdot \nabla C)|_{y=0}=\Da_{\mathrm{w}}\ \chi_{\mathrm{s}}\ [-(1-\wt{B})\wt{C}(x,0,t)+KB]. \label{b bc 1}
\end{equation}
In (\ref{b bc 1}) $\mathbf{n}=(0,-1)$ denotes the outward unit normal vector, $\chi_{\mathrm{s}}$ is the characteristic function defined as
\begin{equation}
\chi_{\mathrm{s}}(\overline{x})=\left\{\begin{array}{ll}
1 & \overline{x}\in [\overline{x}_{\mathrm{min}},\overline{x}_{\mathrm{max}}],\\
0 & \overline{x}\not \in [\overline{x}_{\mathrm{min}},\overline{x}_{\mathrm{max}}],\end{array} \right.  
\end{equation}
and $K=\wt{k}_{\mathrm{d}}/(\wt{k}_{\mathrm{a}}\wt{C}_{\mathrm{u}})$ is the dimensionless equilibrium dissociation rate constant.  Furthermore, since the bound ligand concentration is governed by the kinetics equation \begin{subequations}  
\begin{align}
&\frac{\p B}{\p t}=(1-B)C(x,0,t)-KB\label{react w},\\
&B(x,0)=0, \label{react ic w}
\end{align}\label{kinetics w}\end{subequations}
we can express (\ref{b bc 1}) as 
\begin{equation}
\frac{\p C}{\p y}(x,0,t)=\Da_{\mathrm{w}}\frac{\p B}{\p t}. \label{b bc 2}
\end{equation}
The complete partial differential equation system is given by (\ref{C w inside}), (\ref{kinetics w}), and (\ref{b bc 2}).

In (\ref{b bc 1}) and (\ref{b bc 2}), the important dimensionless parameter 
\begin{equation}
\Da_{\mathrm{w}}=\frac{\wt{H}\wt{k}_{\ra}\wt{R}_{\mathrm{t}}}{\wt{D}}=\frac{\wt{k}_{\ra}\wt{R}_{\mathrm{t}}}{\wt{D}/\wt{H}}\label{Daw}
\end{equation}
is the {\it Damk{\"o}hler number}, which is the ratio of reaction  velocity  to diffusion velocity. Note that both the numerator and denominator have dimensions of unit length per unit time.  It is seen in Table \ref{Table: parameter values} that $\Da_{\mathrm{w}}\gg 1$, which implies that reaction velocity is much faster than diffusion velocity.  This is a direct consequence of the fact that there are multiple time and length scales associated with the experiment: ligand molecules must diffuse a distance on the order of millimeters to arrive at the  biochemical gate, and the speed at which this transpires is far slower than the reaction velocity, \textit{i.e.}, the reaction is diffusion limited. 

Using the fact that $\Da_{\mathrm{w}}\gg 1$ reduces (\ref{b bc 2}) to
\begin{equation}
\frac{\p B}{\p t}=0,\label{leading order steady}
\end{equation}
which implies that, to leading-order, $B(x,t)$ is in a steady-state.  Substituting (\ref{leading order steady}) into (\ref{kinetics w}) yields
\begin{equation}
C(\overline{x},0,t)=\frac{KB}{1-B}.
\end{equation}
This reflects the transport-limited nature of the kinetics system under consideration.  To study the reaction dynamics we must examine the diffusion of ligand molecules in the vicinity of the biochemical gate. We introduce boundary layer coordinates,
\begin{equation}
x=\frac{\overline{x}-1/2}{l_s},\qquad y=\frac{\epsilon}{l_s}\overline{y}. \label{x and y}
\end{equation}
In (\ref{x and y})
\begin{equation}
l_s=\frac{\wt{l}_{\mathrm{s}}}{\wt{L}}
\end{equation}
is the ratio of  biochemical gate length $\wt{l}_{\mathrm{s}}$ to the well length $\wt{L}$, and is very small.  Introducing these scalings into (\ref{C diff w})--(\ref{C no flux w}) and (\ref{b bc 2}) yields
\begin{subequations}
\begin{align}
&\frac{\p C}{\p t}=D\left(\frac{\p^2 C}{\p x^2}+\frac{\p^2 C}{\p y^2}\right),\label{C diff}\\
&C(x,y ,0)=0, \label{ic}\\
&C(x,\epsilon/l_{\mathrm{s}},t)=1,\label{t bc 1}\\
&\frac{\p C}{\p x}(-1/(2l_{\mathrm{s}}),y,t)=\frac{\p C}{\p x}(1/(2l_{\mathrm{s}}),y,t)=0,\label{lr bcs}\\
&\frac{\p C}{\p y}(\overline{x},0,t)=\Da\frac{\p B}{\p t}\chi_{\mathrm{s}}.\label{diff flux}
\end{align}
\end{subequations}
Furthermore, the kinetics equation (\ref{kinetics w}) becomes
\begin{subequations}  
\begin{align}
&\frac{\p B}{\p t}=(1-B)C(x,0,t)-KB\label{react},\\
&B(x,0)=0. \label{react ic}
\end{align}\label{kinetics}\end{subequations}

Observe that transitioning to boundary layer coordinates has the effect of rescaling $D_{\mathrm{w}}$ and $\Da_{\mathrm{w}}$.  The parameter
\begin{equation}
D=\frac{\wt{D}}{\wt{l}_{\mathrm{s}}^2\wt{k}_{\mathrm{a}}\wt{C}_{\mathrm{u}}}=\frac{\wt{D}/\wt{l}_{\mathrm{s}}^2}{\wt{k}_{\ra}\wt{C}_{\mathrm{u}}}
\end{equation}
is the dimensionless diffusion coefficient on this length scale, and is the ratio of the diffusive time scale over a region of size $\wt{l}_{\mathrm{s}}^2$ to the forward reaction rate.  From Table \ref{Table: parameter values} it is seen that $D\gg1$, which implies that diffusion within the boundary layer  is much faster than the forward reaction rate. This is not surprising as we are considering picomolar to femtomolar ligand concentrations.  Furthermore
\begin{equation}
\Da=\frac{\wt{k}_{\ra}\wt{R}_{\mathrm{t}}\wt{l}_{\mathrm{s}}}{\wt{D}}=\frac{\wt{k}_{\mathrm{a}}\wt{R}_{\mathrm{t}}}{\wt{D}/\wt{l}_{\mathrm{s}}}\label{Da}
\end{equation}
is the Damk{\"o}hler number associated with these length scales.  Since $\Da$ is an $O(1)$ to $O(10)$ parameter, on these length scales the reaction velocity is the same as or only slightly faster than the diffusion velocity.  Equation (\ref{diff flux}) then implies that reaction balances  diffusion within the boundary layer.

\subsection{Integrodifferential Equation Reduction} \label{Subesction: Integrodifferential Equation Reduction}

Since $D\gg 1$, we neglect the left hand side of (\ref{C diff}) which reduces this equation to
\begin{equation}
\nabla^2 C=0. \label{quasi-steady C}
\end{equation}
Physically, equation (\ref{quasi-steady C}) implies that near the surface $C$ is in a quasi-steady-state and change in the unbound concentration is driven by the surface-reaction (\ref{diff flux}).  Furthermore, since $l_{\mathrm{s}}\ll 1$ we are not concerned with satisfying the no-flux conditions (\ref{lr bcs}) and take our domain to be the infinite strip $\mathbb{R}\times[0,\epsilon/l_{\mathrm{s}}]$.  This idealization  is physically motivated and justified by the fact that the biochemical gate occupies a very narrow portion of the well surface, so the walls of the well will not appreciably affect ligand binding.

\begin{figure}[tbhp!]
\begin{center}
\begin{tikzpicture}

\node [right] at (4,0) {$\mathrm{Re}\ \omega$};
\node [above] at (0,4) {$\mathrm{Im}\ \omega$};
\node [below] at (1/2+.1,0) {$\rho_n$};
\node [below] at (-1/2-.2,0) {$-\rho_n$};
\node [below] at (3,0) {$R_n$};
\node [below] at (-3,0) {$-R_n$};
\node [above left] at (0,3/4) {$C_2^{(n)}$};
\node [above left] at (0,3) {$C_4^{(n)}$};
\node [above] at (-2,0) {$C_1^{(n)}$};
\node [above] at (2,0) {$C_3^{(n)}$};

\draw[<->](4,0)--(-4,0);
\draw[->](0,0)--(0,4);

\begin{scope}[very thick,decoration={
    markings,
    mark=at position .5 with {\arrow{>}}}
    ]
    \draw [postaction={decorate}] (-3,0)--(-3/4,0);
	\draw [postaction={decorate}] (3/4,0)--(3,0);   
\end{scope}

\begin{scope}[very thick,decoration={
    markings,
    mark=at position 0.25 with {\arrow{<}},
    mark=at position 0.75 with {\arrow{<}}}
    ]

\draw [domain=0:180,postaction={decorate}] plot ({3/4*cos(\x)},{3/4*sin(\x)});
\end{scope}

\begin{scope}[very thick,decoration={
    markings,
    mark=at position 0.25 with {\arrow{>}},
    mark=at position 0.75 with {\arrow{>}}}
    ]
\draw [domain=0:180,postaction={decorate}] plot ({3*cos(\x)},{3*sin(\x)});
\end{scope}

\end{tikzpicture}
\end{center}
\caption{The contour used to calculate (\ref{inv trans}) when $x<0$.}
\label{Figure: contour int 1}
\end{figure}
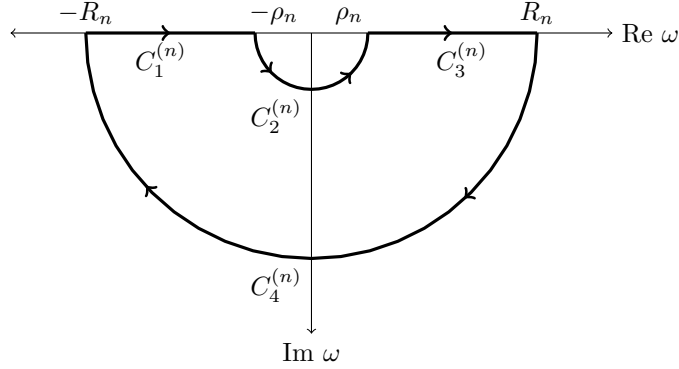

To solve the resulting set of PDEs we seek solutions of the form
\begin{equation}
C(x,y,t)=1+C_{\mathrm{b}}(x,y,t), \label{C ansatz}
\end{equation} 
where $C_{\mathrm{b}}$ satisfies\begin{subequations}
\begin{align}
&\nabla^2C_{\mathrm{b}}=0,\label{Laplace Cb}\\
&C_{\mathrm{b}}(x,\epsilon/l_{\mathrm{s}},t)=0,\label{Cb t bc}\\
&\frac{\p C_{\mathrm{b}}}{\p y}(x,0,t)=\Da\frac{\p B}{\p t}\chi_{\mathrm{s}},\label{Cb b bc}
\end{align}\label{Cb sys}%
\end{subequations}
for $(x,y)\in \mathbb{R}\times[0,\epsilon/l_{\mathrm{s}}]$.
To solve (\ref{Cb sys}) we introduce a Fourier transform in $x$, defining the Fourier transform as\begin{subequations}
\begin{equation}
(\mathcal{F}u)(\omega):=\hat{u}(\omega)=\int_{-\infty}^{\infty}\!u(x)\mathrm{e}^{i\omega x}\ \mathrm{d}x,
\end{equation}
so that the inverse Fourier Transform is given by
\begin{equation}
(\mathcal{F}^{-1}\hat{u})(x)=u(x)=\frac{1}{2\pi}\int_{-\infty}^{\infty}\!\hat{u}(\omega)\mathrm{e}^{-i\omega x}\ \mathrm{d}x.
\end{equation}
\end{subequations}
Applying a Fourier transform to (\ref{Cb sys}) and solving the resulting equations in the frequency domain gives
\begin{equation}
\wh{C}_{\mathrm{b}}(\omega,y,t)=-\frac{\Da\ \sinh((\epsilon l_s^{-1}-y)\omega)}{\omega\cosh(\epsilon l_{s}^{-1}\omega)}\ \frac{\p \wh{B}}{\p t}(\omega,t)\star\left( \frac{\sin(\omega/2)}{\omega/2}\right) \label{f trans sol 2},
\end{equation}
where the convolution product $\star$ has been defined so that
\begin{equation}
\frac{\p \wh{B}}{\p t}(\omega, t)\star\left( \frac{\sin(\omega/2)}{\omega/2}\right)=\int_{-\infty}^{\infty}\frac{\p \wh{B}}{\p t}(\omega-\nu,t)\frac{\sin(\nu/2)}{\nu/2}\ \mathrm{d}\nu.
\end{equation}
However, in order to study the dynamics of interest a closed-form of $C(x,y,t)$ on the surface $y=0$ is required. This is 
aquired by applying the convolution theorem after calculating
\begin{equation}f(x)=\frac{1}{2\pi}\int_{-\infty}^{\infty}\!\frac{\tanh(a\omega)}{\omega}\rme^{-i\omega x}\ \mathrm{d}\omega. \label{inv trans}
\end{equation}
Observe that when $x=0$ the integrand decays at a rate of $1/\omega$ as $\omega\to \pm\infty$.  Thus the integrand of (\ref{inv trans}) is not integrable when $x=0$, and $f$ is singular at the origin.  The evaluation of (\ref{inv trans}) may then be separated into two cases: when $x>0$ and when $x<0$.  We consider the latter by constructing a sequence of contours in the complex plane in the manner depicted in Figure \ref{Figure: contour int 1}. To fix notation we let $C^{(n)}=\sum C_j^{(n)}$.  

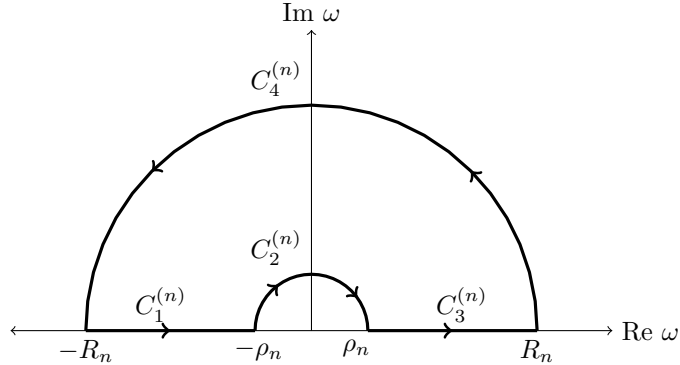
\begin{figure}[tbhp!]
\begin{center}
\begin{tikzpicture}

\node [right] at (4,0) {$\mathrm{Re}\ \omega$};
\node [below] at (0,-4) {$\mathrm{Im}\ \omega$};
\node [above] at (1/2,0) {$\rho_n$};
\node [above] at (-1/2,0) {$-\rho_n$};
\node [above] at (3,0) {$R_n$};
\node [above] at (-3,0) {$-R_n$};
\node [below left] at (0,-3/4) {$C_2^{(n)}$};
\node [below left] at (0,-3) {$C_4^{(n)}$};
\node [below] at (-2,0) {$C_1^{(n)}$};
\node [below] at (2,0) {$C_3^{(n)}$};

\draw[<->](4,0)--(-4,0);
\draw[->](0,0)--(0,-4);

\begin{scope}[very thick,decoration={
    markings,
    mark=at position .5 with {\arrow{>}}}
    ]
    \draw [postaction={decorate}] (-3,0)--(-3/4,0);
	\draw [postaction={decorate}] (3/4,0)--(3,0);   
\end{scope}

\begin{scope}[very thick,decoration={
    markings,
    mark=at position 0.25 with {\arrow{>}},
    mark=at position 0.75 with {\arrow{>}}}
    ]

\draw [domain=0:180,postaction={decorate}] plot ({-3/4*cos(\x)},{-3/4*sin(\x)});
\end{scope}

\begin{scope}[very thick,decoration={
    markings,
    mark=at position 0.25 with {\arrow{<}},
    mark=at position 0.75 with {\arrow{<}}}
    ]
\draw [domain=0:180,postaction={decorate}] plot ({-3*cos(\x)},{-3*sin(\x)});
\end{scope}

\end{tikzpicture}
\end{center}
\caption{The contour used to calculate (\ref{inv trans}) when $x<0$.}
\label{Figure: contour int 2}
\end{figure}

 The hyperbolic tangent function has countably infinite singularities along the imaginary axis, so the path of integration cannot intersect any of these singularities. The singularities will occur when $\omega=0$ or
\begin{equation}
\omega=\frac{\pi i(2n+1)}{2 a}.\label{im singularities}
\end{equation}
Note the contour depicted in Figure \ref{Figure: contour int 1} does not pass through the singularity at the origin; in fact, since
\begin{equation}
\lim_{\omega\to 0}\tanh(a\omega)=0
\end{equation}
this singularity would not have contributed to (\ref{inv trans}) if we  had placed the semi-circle of radius $\rho_n$ in the lower half-plane.  
  Thus taking the radii of our semi-circles to be 
\begin{subequations}
\begin{align}
&\rho_n=\pi/((n+2)a),\\
&R_n=\pi n/a,
\end{align}\label{semi-circle radii}%
\end{subequations}  
the path of integration will never intersect any of the singularities and Cauchy's Residue Theorem may be applied:
\begin{align}
\oint_{C^{(n)}}\!\frac{\tanh(a\omega)}{\omega}\rme^{-i\omega x}\ \mathrm{d}\omega=2\pi i\sum_{k=0}^{n-1}I(C^{(n)},a_k)\ \mathrm{Res}\left(\frac{\tanh(a\omega)}{\omega}\rme^{-i\omega x}; \alpha_n\right).\label{res app}
\end{align}
Calculating residues and letting $n$ approach infinity gives
\begin{equation}
\lim_{n\to \infty}\oint_{C^{(n)}}\!\frac{\tanh(a\omega)}{\omega}\rme^{-i\omega x}\ \mathrm{d}\omega=4\sum_{k=0}^\infty\frac{ \rme^{(2k+1)\pi x/(2 a)}}{(2k+1)}.\label{res sum}
\end{equation}
On the other hand,
\begin{equation}
\lim_{n\to\infty}\oint_{C^{(n)}}\! \frac{\tanh(a\omega)}{\omega}\ \mathrm{e}^{-i\omega x}\ \mathrm{d}\omega=\lim_{n\to\infty}\sum_{j=1}^4\oint_{C_j^{(n)}}\! \frac{\tanh(a\omega)}{\omega}\ \mathrm{e}^{-i\omega x}\ \mathrm{d}\omega.\label{contour sum}
\end{equation}
One may show that the integral along $C_2^{(n)}$ vanishes as $n\to\infty$, and by using the fact that $x<0$ one may similarly show that the integral along the far contour $C_2^{(4)}$ vanishes.  From these facts and the Maclaurin series for $\tanh^{-1}(x)$ it follows that
\begin{equation}
f(x)=\frac{2}{\pi}\tanh^{-1}(\rme^{\pi l_s x/(2\epsilon)})\label{f inv x<0}
\end{equation}
 when $x<0$.  To evaluate (\ref{inv trans}) when $x>0$ one may extend this integral to the complex plane by using the reflection of the contour depicted in Figure \ref{Figure: contour int 1} about the real axis, shown in Figure \ref{Figure: contour int 2}, and use analogous arguments to show
\begin{equation}
f(x)=\frac{2}{\pi}\tanh^{-1}(\rme^{-\pi l_s x/(2\epsilon)})\label{f inv x>0}
\end{equation}
when $x>0$.  

In summary the integral (\ref{inv trans}) is singular at the origin, given by (\ref{f inv x<0}) when $x<0$, and  (\ref{f inv x>0}) when $x>0$.  Putting these three observations together leads to the conclusion that 
\begin{equation}
f(x)=\tanh^{-1}(\rme^{-\pi l_s|x|/(2\epsilon)}).\label{f inv}
\end{equation}
Thus applying the convolution theorem to (\ref{f trans sol 2}) evaluated at $y=0$ and substituting the resulting expression into (\ref{C ansatz}) gives:
\begin{equation}
C(x,0,t)=1-\frac{2\ \Da}{\pi}\int_{-1/2}^{1/2}\!\tanh^{-1}(\rme^{-\pi l_s|x-\nu|/(2\epsilon)})\frac{\p B}{\p t}(\nu,t)\ \mathrm{d}\nu.\label{surf C}
\end{equation}
Hence, the bound ligand concentration is governed by the IDE\begin{subequations}
\begin{align}
&\frac{\p B}{\p t}=(1-B)\left(1-\frac{2\ \Da}{\pi}\int_{-1/2}^{1/2}\!\tanh^{-1}(\rme^{-\pi l_s|x-\nu|/(2\epsilon)})\frac{\p B}{\p t}(\nu,t)\ \mathrm{d}\nu\right)-KB,\label{ide}\\%
&B(x,0)=0.\label{ide ic}%
\end{align}\label{ide prob}%
\end{subequations}%
In (\ref{surf C}) the term 1 represents the uniform injection concentration and the convolution integral represents depletion of unbound ligand at the surface due to reaction.  As we shall see in Section \ref{Section: Results and Discussion} the non-local nature of the convolution (\ref{surf C}) reflects the probabilistic nature of diffusion in the boundary layer near the surface, and the finite limits of integration encode the reflective boundary conditions to the left and right of the biochemical gate.  However, we first turn our attention to finding a numerical approximation to the solution of (\ref{ide prob}).

\section{Numerical Method} \label{Section: Numerical Method}

\subsection{Method of Lines Approximation} \label{Subsection: Method of Lines Approximation}

To discretize (\ref{ide}) we choose $N$ equally-spaced discretization nodes $x_i$ and partition $[-1/2,1/2]$ into $N$ distinct subintervals of length $\dx=1/N$:
\begin{equation}
\left[-\frac{1}{2},\frac{1}{2}\right]=\bigcup_{i=1}^N \left[x_i-\frac{\dx}{2},x_i+\frac{\dx}{2}\right],
\end{equation} 
where $-1/2=x_1-\dx/2$ and $1/2=x_N+\dx/2$.  Then  an approximation to (\ref{ide prob}) is found by applying the method of lines
\begin{equation}
B(x,t)\approx\sum_{i=1}^nh_i(t)\phi_i(x) \label{mol exp}
\end{equation}
where the functions $h_i(t)$ are to be determined and subject to the initial condition $h_i(0)=0$, and the functions $\phi_i(x)$ are locally defined piece-wise linear hat functions
\begin{equation}
\phi_i(x)=\left\{\begin{array}{ll}
\displaystyle\frac{2}{\dx}[x-(x_i-\dx/2)] & \mathrm{if }\ x\in[x_i-\dx/2,x_i),\vspace{1em}\\
\displaystyle\frac{2}{\dx}[(x_i+\dx/2)-x] & \mathrm{if }\ x\in[x_i,x_i+\dx/2],\vspace{1em}\\
0 & \mathrm{else} .
\end{array}\right. \label{phi_i}
\end{equation}
Substituting (\ref{mol exp}) into (\ref{ide}) and evaluating each side of the resulting equation at $x=x_j$ yields
\begin{align}
\begin{aligned}
h_j'(t)=&\left(1-h_j(t)\right)\left(1-\sum_{i=1}^N\frac{2\ \Da\ h_i'(t)}{\pi}\int_{-1/2}^{1/2}\tanh^{-1}(\rme^{-|x_j-\nu|\pi l_{\mathrm{s}}/(2\epsilon)})\phi_i(\nu)\ \mathrm{d}\nu\right)\\&\quad-Kh_j(t), 
\end{aligned}\label{mol ode 3}
\end{align}
for $j=1,\ \ldots,\ N$.  The solution of this nonlinear set of ODEs  determines the time-dependent functions $h_j(t)$, however solving this system requires computing
\begin{equation}
\int_{-1/2}^{1/2}\tanh^{-1}(\rme^{-|x_j-\nu|\pi l_{\mathrm{s}}/(2\epsilon)})\phi_i(\nu)\ \mathrm{d}\nu.\label{inv tanh int}
\end{equation}
Since $\tanh^{-1}(\rme^{-|x_j-\nu|\pi l_{\mathrm{s}}/(2\epsilon)})$ exhibits logarithmic singularity at $\nu=x_j$, computing (\ref{inv tanh int}) using a quadrature rule requires great care, although (\ref{inv tanh int}) may be evaluated exactly.  This is done by decomposing the basis functions (\ref{phi_i}) into their left and right parts:
\begin{equation}
\phi_{i,l}(x)=\left\{\begin{array}{ll}
\displaystyle\frac{2}{\dx}[x-(x_i-\dx/2)] & \mathrm{if }\ x\in[x_i-\dx/2,x_i),\vspace{1em}\\
0& \mathrm{else}, 
\end{array}\right. 
\end{equation}
and
\begin{equation}
\phi_{i,r}(x)=\left\{\begin{array}{ll}
\displaystyle\frac{2}{\dx}[(x+x_i)-\dx/2] & \mathrm{if }\ x\in[x_i,x_i+\dx/2],\vspace{1em}\\
0& \mathrm{else}.
\end{array}\right. 
\end{equation}
Having decomposed the basis functions into their left and right parts (\ref{inv tanh int}) can be written as 
\begin{align}\begin{aligned}
\int_{-1/2}^{1/2}\tanh^{-1}(\rme^{-|x_j-\nu|\pi l_{\mathrm{s}}/(2\epsilon)})\phi_i(\nu)\ \mathrm{d}\nu=&\int_{-1/2}^{1/2}\tanh^{-1}(\rme^{-|x_j-\nu|\pi l_{\mathrm{s}}/(2\epsilon)})\phi_{i,l}(\nu)\ \mathrm{d}\nu\\ &\ +\int_{-1/2}^{1/2}\tanh^{-1}(\rme^{-|x_j-\nu|\pi l_{\mathrm{s}}/(2\epsilon)})\phi_{i,r}(\nu)\ \mathrm{d}\nu.\end{aligned}\label{inv tanh int sum}
\end{align}
Since the two terms on the right hand side are related through a change of variables, it is sufficient to calculate
\begin{equation}
\int_{-1/2}^{1/2}\tanh^{-1}(\rme^{-|x_j-\nu|\pi l_{\mathrm{s}}/(2\epsilon)})\phi_{i,l}(\nu)\ \mathrm{d}\nu. \label{phi l int}
\end{equation}
After changing variables, one may use the definition of $\tanh^{-1}(\cdot)$ and expand the integrand in terms of its Mclaurin series to find that it is a telescoping sum:
\begin{align}
&\int_{-1/2}^{1/2}\tanh^{-1}(\rme^{-|x_j-\nu|\pi l_{\mathrm{s}}/(2\epsilon)})\phi_{i,l}(\nu)\ \mathrm{d}\nu\\ &\quad =\sum_{n=0}^{\infty}\frac{2}{\Delta x}\int_0^{\dx/2}\frac{\rme^{-|w-x_j+x_i-\Delta x/2|(2n+1)/(2\epsilon)}}{2n+1}w\ \mathrm{d}w. \label{telescope}
\end{align}
In writing (\ref{telescope}) we have formally exchanged the limit operations.  Observe that the absolute value prevents one from integrating by parts directly; however, by using the fact that the discretization nodes are equally spaced one can show the computation may be partitioned in two distinct cases: when $x_j\ge x_i$ and $x_j<x_i$.  Since the computation is analogous in each case we concern ourselves only with the former.  Thus taking $x_j\ge x_i$ and integrating the right hand side of (\ref{telescope}) by parts shows that (\ref{phi l int}) is equal to
\begin{align}
\begin{aligned}
&\sum_{n=0}^\infty\left(\frac{2}{\dx}\right)\bigg(\frac{\dx \epsilon\  \rme^{-(x_j-x_i)(2n+1)\pi l_{\mathrm{s}}/(2\epsilon)}}{(2n+1)^2\pi l_{\mathrm{s}}}-\frac{4\epsilon^2\ \rme^{-(x_j-x_i)(2n+1)\pi l_{\mathrm{s}}/(2\epsilon)}}{(2n+1)^3\pi^2 l_{\mathrm{s}}^2}\\&\quad+\frac{4\epsilon^2\  \rme^{-[\dx/2+(x_j-x_i)](2n+1)\pi l_{\mathrm{s}}/(2\epsilon)}}{(2n+1)^3\pi^2 l_{\mathrm{s}}^2}\bigg).\label{horrible series}\end{aligned}
\end{align}
To sum the series (\ref{horrible series}), we observe that one can use the definition of the polylogarithm of order $s$
\begin{equation}
\mathrm{Li}_{\mathrm{s}}(z)=\sum_{k=1}^\infty\frac{z^k}{k^s}
\end{equation}
to show
\begin{equation}
\sum_{n=0}^\infty\frac{z^{2n+1}}{(2n+1)^s}=\mathrm{Li}_{\mathrm{s}}(z)-\frac{1}{2^{\mathrm{s}}}\mathrm{Li}_{s}(z^2).
\end{equation}
Hence when $x_j\ge x_i$
\begin{align}\begin{aligned}
&\int_{-1/2}^{1/2}\tanh^{-1}(\rme^{-|x_j-\nu|\pi l_{\mathrm{s}}/(2\epsilon)})\phi_{i,l}(\nu)\ \mathrm{d}\nu\\ &\quad =\left(\frac{2}{\Delta x}\right)\bigg[\frac{\dx \epsilon}{\pi l_{\mathrm{s}}}\left(\mathrm{Li}_2(\rme^{-(x_j-x_i)\pi l_{\mathrm{s}}/(2\epsilon)})-\mathrm{Li}_2(\rme^{-(x_j-x_i)\pi l_{\mathrm{s}}/(\epsilon)})/4\right)\\&\qquad-\frac{4 \epsilon^2}{\pi^2 l_{\mathrm{s}}^2}\left(\mathrm{Li}_3(\rme^{-(x_j-x_i)\pi l_{\mathrm{s}}/(2\epsilon)})-\mathrm{Li}_3(\rme^{-(x_j-x_i)\pi l_{\mathrm{s}}/(\epsilon)})/8\right)\\&\qquad +\frac{4 \epsilon^2}{\pi^2 l_{\mathrm{s}}^2}\left(\mathrm{Li}_3(\rme^{-[\dx/2+(x_j-x_i)]\pi l_{\mathrm{s}}/(2\epsilon)})-\mathrm{Li}_3(\rme^{-[\dx/2+(x_j-x_i)]\pi l_{\mathrm{s}}/(\epsilon)})/8\right)\bigg].\end{aligned} \label{phi l int val}
\end{align}
The form of (\ref{phi l int val}) when $x_j<x_i$ is similar.  With the exact value of (\ref{inv tanh int sum}), the nonlinear set of ODEs (\ref{mol ode 3}) may be integrated with one's favorite linear multistage or multistep formula.

\subsection{Convergence} \label{Subsection: Convergence}

\begin{figure}[tbhp!]
  \centering
  \includegraphics[width=0.5\textwidth]{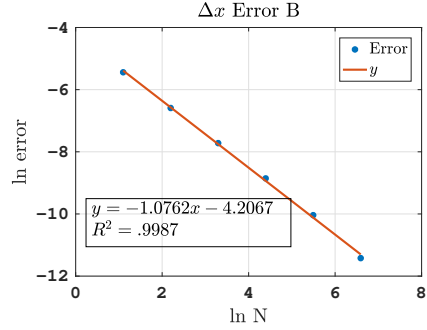}
  \caption{The values of (\ref{error}) for $i=1,\ \ldots,\ 6$ depicted together with the line $y=-1.0762x-4.2067$, which was fit to the values of (\ref{error}) with an $R^2$ coefficient of $R^2=.9987$.  Parameter values of $\Da=66.42$, $K=1$, $l_{\mathrm{s}}=10^{-3}$, and $\epsilon=1$ were used.}\label{Figure: conv fig}
\end{figure}

Convergence of the numerical method outlined in the previous subsection was measured by first computing a reference solution $B_{\mathrm{ref}}(x,t)$ on a mesh with $N=3^7=2187$ spatial discretization nodes; this was done by integrating   (\ref{mol ode 3}) from $t=0$ to $t=150$ using an adaptive linear multistage formula.  Then  solutions $B_i(x,t)$ were computed on meshes with $N=3^i$ nodes and  convergence was measured by calculating
\begin{align}
&||\ ||B_{\mathrm{ref}}(x,t)-B_i(x,t)||_{2,\ x}||_{\infty,\ t}\label{error}
\end{align}
for $i=1,\ \ldots, 6$.  In (\ref{error}) $|| \cdot||_{2,\ x}$ denotes $l_2$ norm in $x$ and $||\cdot||_{\infty,\ t}$ denotes the infinity norm in $t$.  A logarithmic plot of these values is depicted in Figure \ref{Figure: conv fig}.  Despite the logarithmic singularity in (\ref{ide}), the evidence in Figure \ref{Figure: conv fig} strongly suggests that our method of lines approximation to (\ref{ide prob}) achieves first-order convergence.  Although it is of interest to derive analytic error estimates for our approximation, the nonlinearity in (\ref{ide}) precludes analysis.

\section{Results and Discussion}\label{Section: Results and Discussion}

\begin{figure}[tbhp!]
  \centering
  \subfloat[Space-time curve of $B(x,t)$ for $t$ in the interval ${[}0,.1{]}$ ]{\includegraphics[width=0.45\textwidth]{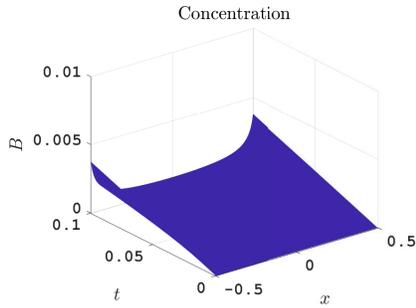}\label{Figure: Bt1}}\qquad
   \subfloat[Space-time curve of $B(x,t)$ for $t$ in the interval ${[}0,10{]} $.]{\includegraphics[width=0.45\textwidth]{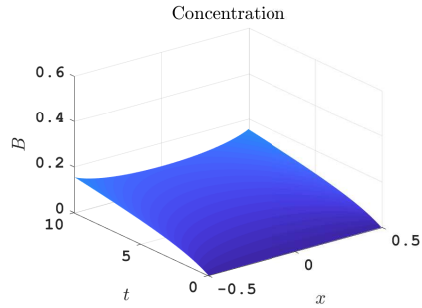}\label{Figure: Bt10}}\hfill
    \subfloat[Space-time curve of $B(x,t)$ for $t$ in the interval ${[}0,50{]}$. ]{\includegraphics[width=0.45\textwidth]{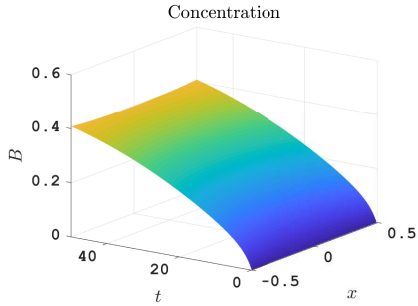}\label{Figure: Bt50}}\qquad
  \subfloat[Space-time curve of $B(x,t)$ for $t$ in the interval ${[}0,150{]}$.]{\includegraphics[width=0.45\textwidth]{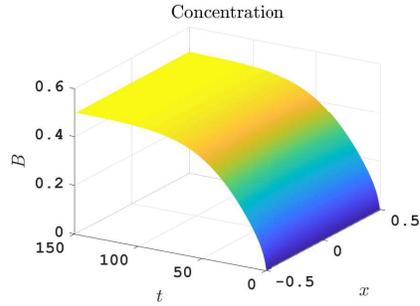}\label{Figure: Bt150}}
  \caption{Method of lines approximation to the solution of (\ref{ide prob}) during different time intervals.  Parameter values of $\Da=66.42$, $K=1$, $l_{\mathrm{s}}=10^{-3}$, and $\epsilon=2/5$ were used.}\label{Figure: B evo}\hfill
\end{figure}

The results of our numerical simulations are depicted in Figure \ref{Figure: B evo}.  Upon inspection one immediately notices the presence of a depletion region in the center of the biochemical gate for small $t$. As time progresses, the depletion regions narrows and becomes more shallow as the rate of bound ligand production near the boundary decreases.   The bound ligand concentration continues to become more spatially uniform  until a chemical equilibrium is achieved, resulting in  a balance between association and dissociation.

\begin{figure}[tbhp!]
    \centering
    \includegraphics[width=0.5\textwidth]{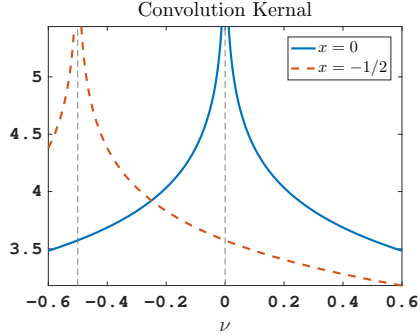}
    \caption{The convolution kernal $\tanh^{-1}(\rme^{-|x-\nu|\pi l_{\mathrm{s}}/(2\epsilon)})$ centered at $x=0$ (solid line), and at $x=-1/2$ (dotted line).  Parameter values of $l_{\mathrm{s}}=10^{-3}$ and $\epsilon=2/5$ were used.}
    \label{Figure: convolution kernal}
\end{figure}

Mathematically, the depletion region results from the singular convolution kernal \begin{equation}
\tanh^{-1}(\rme^{-|(x-\nu)|\pi l_{\mathrm{s}}/(2\epsilon)})
\end{equation} 
and the finite limits of integration.  In Figure \ref{Figure: convolution kernal} the convolution kernel has been depicted,  centered at both $x=0$ and $x=-1/2$.  When the convolution  kernel is centered at $x=0$ it acts as a two-sided influence function.  The singularity at $x=0$ reflects the high  likelihood that a ligand molecule directly above the origin will diffuse to the surface and bind with an available receptor site there; however, in the unstirred layer ligand molecules diffusing into the surface bind with  neighboring receptor sites.  Figure \ref{Figure: convolution kernal} reveals the likelihood of binding with  a neighboring receptor site decays with the distance away from the source, although it is never zero since  $\tanh^{-1}({\rme^{-|x-\nu|\pi l_{\mathrm{s}}}/(2\epsilon)})$ is supported everywhere on the real line, and in particular everywhere on $[-1/2,1/2]$.  Conversely, when the kernel is centered at $x=-1/2$ Figure \ref{Figure: convolution kernal} shows that it acts as a one-sided influence function. The finite limits of integration in (\ref{ide}) imply that the convolution kernel influences the bound ligand concentration the most at $x=-1/2$, and has a monotonically decreasing influence progressing from $x=-1/2$ to $x=1/2$.  Thus the finite limits of integration encode the reflective boundary conditions.  To the right of $x=-1/2$ ligand molecules spread out and diffuse into the surface, while to the left they are merely reflected. 

\begin{figure}[tbhp!] 
  \centering
  \subfloat[The average concentration has been depicted for $\wt{k}_{\mathrm{a}}=10^{11},\ 5\times10^{11}$, and $10^{12}\ \mathrm{cm}^3/(\mathrm{mol}\cdot \mathrm{s})$.  This corresponded to $\Da=6.64,\ 33.21,$ and $\ 66.42$; and $K=1.67,\ 0.33,$ and  $0.17$.  In addition parameter values of $l_{\mathrm{s}}=10^{-3}$, and $\epsilon=2/5$ were used. ]{\includegraphics[width=0.46\textwidth]{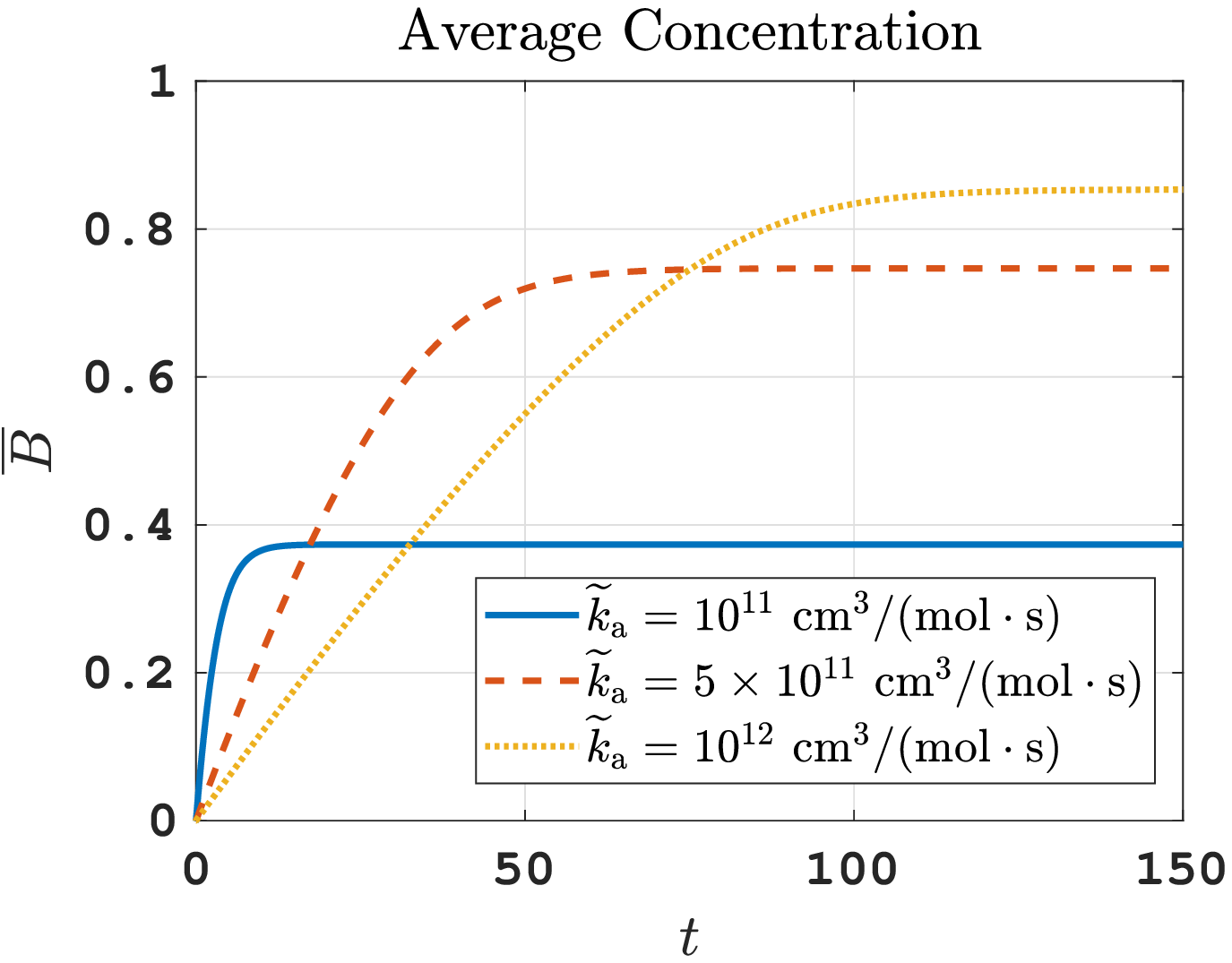}\label{Figure: Bbar ka}}\qquad
  \subfloat[The average concentration has been depicted for $\wt{C}_{\mathrm{u}}=10^{-17},\ 5\times10^{-17}$, and $10^{-16}\ \mathrm{mol}/\mathrm{cm}^3$.  This corresponded to $K=10,\ 2,$ and $0.2$.  In addition parameter values of $\Da=6.6420,\ l_{\mathrm{s}}=10^{-3}$, and $\epsilon=2/5$ were used. ]{\includegraphics[width=0.46\textwidth]{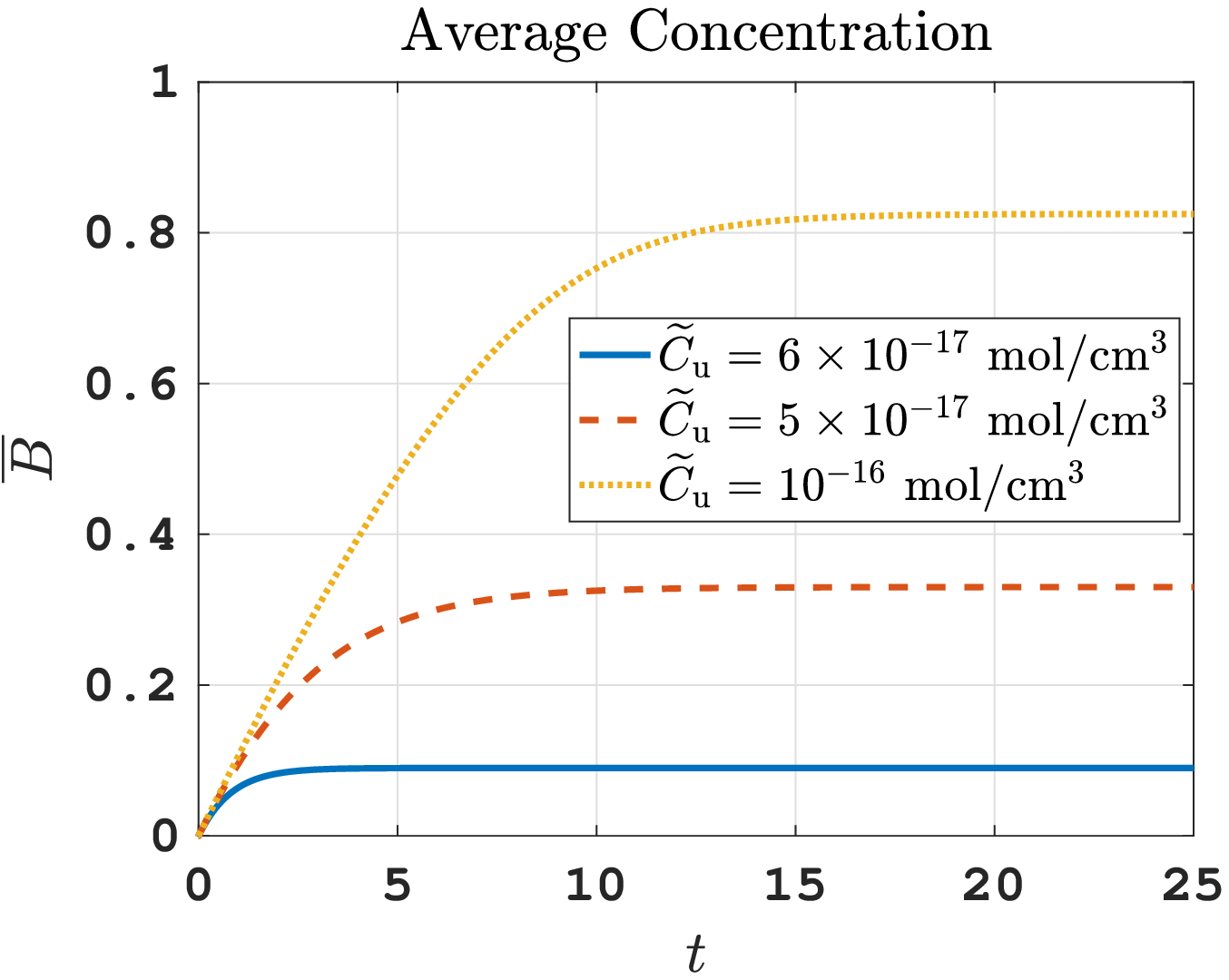}\label{Figure: Bbar Cu}}
   \caption{The average concentration (\ref{average concentration}) for different values of $\wt{k}_{\mathrm{a}}$ and $\wt{C}_{\mathrm{u}}$.}
    \label{Figure: Bbar}
\end{figure}

The average concentration across the biochemical gate
\begin{equation}
\overline{B}(t)=\int_{-1/2}^{1/2}B(x,t)\ \mathrm{d}x\label{average concentration}
\end{equation}
is shown in Figure \ref{Figure: Bbar ka} for three values of $\wt{k}_{\mathrm{a}}$. This quantity is proportional to the electrostatic potential applied to the biochemical gate, and thereby the electric current across the semi-conducting channel, allowing direct comparison to measurements.  Increasing the association rate constant results in a larger Damk{\"o}hler number. This reflects the enhanced rate of reaction relative to transport, and corresponds to wider and deeper depletion regions that impede current flow near the boundaries of the biochemical gate before the rest of the semiconductor channel. This is a remarkable result that is not directly observable experimentally, and provides physical insight into the origin of the signal measured with a FET. Finally,  the transient phase of the signal grows with the association rate constant, owing to both decreasing the equilibrium dissociation rate constant and increasing the rate of reaction relative to diffusion. 

Increasing the ligand concentration $\wt{C}_{\mathrm{u}}$ increases  the average concentration at the biochemical gate, resulting in higher   FET conductance.  From Figure \ref{Figure: Bbar Cu} it is seen that the  equilibrium value of $\overline{B}$ increases with the ligand concentration $\wt{C}_{\mathrm{u}}$, as expected since the ligand concentration and equilibrium dissociation rate constant are inversely proportional.  These considerations are clearly of fundamental importance for parameter estimation.  



\section{Conclusions}\label{Section: Conclusions}

The ability to tailor therapies to individuals or specific subsets of a population could transform medicine.  However,  widespread use of personalized therapeutics has yet to be adopted due to our inability to quickly and routinely measure biomarkers.  Not only do FETs exhibit high charge sensitivity and provide direct signal transduction, they also provide label-free measurements at physiological concentrations.  As such, FETs are  an incredibly promising tool for biomarker measurement.  Although an accurate dynamical model for receptor-ligand dynamics is necessary for maximizing the sensitivity of these instruments, all previous modeling efforts have been devoted to the study of steady-state sensor physics.  Herein, a time-dependent model for receptor ligand dynamics has been presented for the first time.

This model takes the form of a diffusion equation, coupled to an equation describing reaction on the biochemical gate.  Analysis of this set of nonlinear equations is complicated by the presence of multiple disparate time and length scales: ligand molecules must diffuse a distance on the centimeters to arrive at the reacting surface, which is on the order of micrometers.  Furthermore, diffusion is a very slow process while the reactions of interest proceed very quickly.  Nevertheless, by using the appropriate characteristic time and length scales one is able to reduce this model to a quasi-steady transport equation for the unbound ligand concentration $C$, coupled to an equation describing the evolution of the bound ligand concentration $B$.  Employing the residue theorem allows one to further reduce this set of equations to a single nonlinear IDE in terms of the reacting species concentration.  Despite the presence of a singular convolution kernel, this equation has been solved to first-order accuracy without the need to resort to specialized quadrature techniques to evaluate (\ref{inv tanh int}).  Results of our numerical simulations reveal the presence of a depletion region in the center of the biochemical gate, which influences the current signal by non-uniformly altering the surface-potential of the semiconductor channel.



In addition to providing a time-dependent model for estimating binding affinities, the present model could be coupled to a model for semiconductor physics to refine theoretical predictions and serve as a basis for sensor optimization.  The latter may be a subject of future investigation.  Additionally, it is of interest to study receptor-ligand dynamics in FETs under a different experimental conditions; \textit{i.e.} a sealed experiment wherein a drop of ligand molecules is injected at an instance of time.  Extending the present model to higher geometries is also of interest.

\section*{Acknowledgements}

 The authors are grateful to Paul Patrone for the many valuable conversations.

\bibliographystyle{siamplain}
\bibliography{research_refs}

\begin{thebibliography}{10}

\bibitem{Armani:2003gx}
D.~K. Armani, T.~J. Kippenberg, S.~M. Spillane, and K.J. Vahala.
\newblock Ultra-high-{Q} toroid microcavity on a chip.
\newblock {\em Letters to Nature}, 421(6926):925, 2003.

\bibitem{Binukumar:2014hl}
B.~B.~K., Y.-L. Zheng, V.~Shukla, N.~D. Amin, P.~Grant, and H.~C. Pant.
\newblock {TFP5}, a peptide derived from {P35}, a {CDK5} neuronal activator,
  rescues cortical neurons from glucose toxicity.
\newblock {\em Journal of Alzheimer's Disease}, 39(4):899--909, 2014.

\bibitem{baumgartner2011analysis}
S.~Baumgartner, M.~Vasicek, A.~Bulyha, N.~Tassotti, and C.~Heitzinger.
\newblock Analysis of field-effect biosensors using self-consistent {3D}
  drift-diffusion and {M}onte-{C}arlo simulations.
\newblock {\em Procedia Engineering}, 25:407--410, 2011.

\bibitem{Boyd:2001tg}
R.~W. Boyd and J.~E. Heebner.
\newblock Sensitive disk resonator photonic biosensor.
\newblock {\em Applied Optics}, 40(31):5742--5747, 2001.

\bibitem{Cardone:2013hh}
A.~Cardone, H.~Pant, and S.~A. Hassan.
\newblock Specific and non-specific protein association in solution:
  computation of solvent effects and prediction of first-encounter modes for
  efficient configurational bias monte carlo simulations.
\newblock {\em The Journal of Physical Chemistry B}, 117(41):12360--12374,
  2013.

\bibitem{Cui:2001iy}
Y.~Cui, Q.~Wei, H.~Park, and C.~M. Lieber.
\newblock Nanowire nanosensors for highly sensitive and selective detection of
  biological and chemical species.
\newblock {\em Science}, 293(5533):1289--1292, 2001.

\bibitem{Dhavan:2001gd}
R.~Dhavan and L.-H. Tsai.
\newblock A decade of cdk5.
\newblock {\em Nature reviews. Molecular Cell Biology}, 2(10):749, 2001.

\bibitem{Drain:2016fs}
P.~K. Drain, L.~Gounder, F.~Sahid, and M.-Y. Moosa.
\newblock Rapid urine {LAM} testing improves diagnosis of expectorated
  smear-negative pulmonary tuberculosis in an {HIV}-endemic region.
\newblock {\em Scientific reports}, 6:19992, 2016.

\bibitem{Fosgerau:2015kt}
K.~Fosgerau and T.~Hoffmann.
\newblock Peptide herapeutics: current status and future directions.
\newblock {\em Drug Discovery Today}, 20(1):122--128, 2015.

\bibitem{heitzinger2010multiscale}
C.~Heitzinger, N.~J. Mauser, and C.~Ringhofer.
\newblock Multiscale modeling of planar and nanowire field-effect biosensors.
\newblock {\em SIAM Journal on Applied Mathematics}, 70(5):1634--1654, 2010.

\bibitem{Henrich:2010ds}
S.~Henrich, O.~Salo-Ahen, B.~Huang, F.~F. Rippmann, G.~Cruciani, and R.~C.
  Wade.
\newblock Computational approaches to identifying and characterizing protein
  binding sites for ligand design.
\newblock {\em Journal of Molecular Recognition}, 23(2):209--219, 2010.

\bibitem{Ilic:2004go}
B.~Ilic, H.~G. Craighead, S.~Krylov, W.~Senaratne, C.~Ober, and P.~Neuzil.
\newblock Attogram detection using nanoelectromechanical oscillators.
\newblock {\em Journal of Applied Physics}, 95(7):3694--3703, 2004.

\bibitem{Johannsmann:2015kq}
D.~Johannsmann and G.~Brenner.
\newblock Frequency shifts of a quartz crystal microbalance calculated with the
  frequency-domain lattice--boltzmann method: application to coupled liquid
  mass.
\newblock {\em Analytical chemistry}, 87(14):7476--7484, 2015.

\bibitem{Knoll:1998df}
W.~Knoll.
\newblock Interfaces and thin films as seen by bound electromagnetic waves.
\newblock {\em Annual Review of Physical Chemistry}, 49(1):569--638, 1998.

\bibitem{landheer2005model}
D.~Landheer, G.~Aers, W.~R. McKinnon, M.~J. Deen, and J.~C. Ranuarez.
\newblock Model for the field effect from layers of biological macromolecules
  on the gates of metal-oxide-semiconductor transistors.
\newblock {\em Journal of Applied Physics}, 98(4):044701, 2005.

\bibitem{Lawn:2016ed}
S.~D. Lawn and A.~Gupta-Wright.
\newblock Detection of lipoarabinomannan ({LAM}) in urine is indicative of
  disseminated tb with renal involvement in patients living with {HIV} and
  advanced immunodeficiency: evidence and implications.
\newblock {\em Transactions of the Royal Society of Tropical Medicine and
  Hygiene}, 110(3):180--185, 2016.

\bibitem{Mohanty:2014wy}
P.~{Mohanty}, Y.~{Chen}, X.~{Wang}, M.~K. {Hong}, C.~L. {Rosenberg}, D.~T.
  {Weaver}, and S.~{Erramilli}.
\newblock {Field Effect Transistor Nanosensor for Breast Cancer Diagnostics}.
\newblock {\em ArXiv e-prints}, 2014.

\bibitem{Naik:2009gu}
A.~K. Naik, M.~S. Hanay, W.~K. Hiebert, X.~L. Feng, and M.~L. Roukes.
\newblock Towards single-molecule nanomechanical mass spectrometry.
\newblock {\em Nature Nanotechnology}, 4(7):445--450, 2009.

\bibitem{Pouthas:2004cj}
F.~Pouthas, C.~Gentil, D.~C{\^o}te, and U.~Bockelmann.
\newblock {DNA} detection on transistor arrays following mutation-specific
  enzymatic amplification.
\newblock {\em Applied Physics Letters}, 84(9):1594--1596, 2004.

\bibitem{Rodahl:1995ck}
M.~Rodahl, F.~H{\"o}{\"o}k, A.~Krozer, P.~Brzezinski, and B.~Kasemo.
\newblock Quartz crystal microbalance setup for frequency and {Q}-factor
  measurements in gaseous and liquid environments.
\newblock {\em Review of Scientific Instruments}, 66(7):3924--3930, 1995.

\bibitem{Su:2015fo}
J.~Su.
\newblock Label-free single exosome detection using frequency-locked
  microtoroid optical resonators.
\newblock {\em ACS Photonics}, 2(9):1241--1245, 2015.

\bibitem{tulzer2013kinetic}
G.~Tulzer, S.~Baumgartner, E.~Brunet, G.~C. Mutinati, S.~Steinhauer,
  A.~K{\"o}ck, P.~E. Barbano, and C.~Heitzinger.
\newblock Kinetic parameter estimation and fluctuation analysis of co at
  sno${}_2$ single nanowires.
\newblock {\em Nanotechnology}, 24(31):315501, 2013.

\bibitem{Walsh:2014kw}
G.~Walsh.
\newblock Biopharmaceutical benchmarks 2014.
\newblock {\em Nature Biotechnology}, 32(10):992--1000, 2014.

\bibitem{Wang:2010bo}
S.~Wang, X.~Shan, U.~Patel, X.~Huang, J.~Lu, J.~Li, and N.~Tao.
\newblock Label-free imaging, detection, and mass measurement of single viruses
  by surface plasmon resonance.
\newblock {\em Proceedings of the National Academy of Sciences},
  107(37):16028--16032, 2010.

\bibitem{Wang:2005ij}
W.~U. Wang, C.~Chen, K.-H. Lin, Y.~Fang, and C.~M. Lieber.
\newblock Label-free detection of small-molecule--protein interactions by using
  nanowire nanosensors.
\newblock {\em Proceedings of the National Academy of Sciences of the United
  States of America}, 102(9):3208--3212, 2005.

\end{thebibliography}

\end{document}